\documentclass[11pt]{amsart} %
\usepackage{amsmath,amssymb, amstext, enumerate}
\usepackage{fullpage}
\usepackage{graphicx}
\usepackage{graphics}
\usepackage{psfrag}

\newtheorem{theorem}{Theorem}[section]

\newtheorem{proposition}[theorem]{Proposition}

\newtheorem{lemma}[theorem]{Lemma}

\numberwithin{equation}{section}

\def\Lim{{\rm Lim}\,}

\def\b0{{\bf 0}}

\def\x{{\bf x}}
\def\y{{\bf y}}

\def\cN{{\mathcal N}}
\def\cV{{\mathcal V}}

\def\diag{{\rm diag}\,}
\def\tr{{\rm tr}\,}

\def\ker{{\rm ker}\,}

\def\({\left(}
\def\){\right)}
\def\[{\left[}
\def\]{\right]}

\def\tr{{\rm tr}}

\begin{document}
\openup 1.3\jot

\title[quadratic operator]{Factoring a quadratic operator as a product \\ of two positive contractions}

\author{Chi-Kwong Li and Ming-Cheng Tsai}

\address[Li]{Department of Mathematics, College of William \& Mary, Williamsburg, VA 23187, USA.}
\email{ckli@math.wm.edu}

\address[Tsai]{Department of Applied Mathematics, National Sun Yat-sen University, Kaohsiung 80424, Taiwan.}
\email{mctsai2@gmail.com}


\keywords{quadratic operator, positive contraction, spectral theorem}

\subjclass{47A60, 47A68, 47A63}
\date{}
\maketitle
\begin{abstract}
Let $T$ be a quadratic operator on a complex Hilbert space $H$. We show that $T$ can be written as a product of two positive contractions if and only if $T$ is of the form
$$aI \oplus bI \oplus\begin{pmatrix} aI & P \cr 0 & bI \cr \end{pmatrix} \quad \text{on} \quad  H_1\oplus H_2\oplus (H_3\oplus H_3)$$
for some $a, b\in [0,1]$ and strictly positive operator $P$ with $\|P\| \le  |\sqrt{a} - \sqrt{b}|\sqrt{(1-a)(1-b)}.$ Also, we give a necessary condition for a bounded linear operator $T$ with operator matrix
$\begin{pmatrix} T_1 & T_3\\ 0 & T_2\cr\end{pmatrix}$ on $H\oplus K$ that can be written as a product of two positive contractions.
\end{abstract}

\maketitle

\section{Introduction}

There has been considerable interest in studying factorization of bounded linear operators (see \cite{Ba1, Ba2, Ba3, Ba4, Wu89}). For example, a $2\times 2$ matrix $C$ can be written as a product of two orthogonal projections if and only if $C$ is the identity operator or $C$ is unitarily similar to $\begin{pmatrix}a & \sqrt{a(1-a)} \cr 0 & 0\cr\end{pmatrix}$ for some $a\in [0,1]$. For more results about products of orthogonal projections, one may consult \cite{Am, Bott, Corach, Ha}. Note that one can write an $n\times n$ matrix $C$ as a product of two positive operators exactly when $C$ is similar to a positive operator (see \cite[Theorem 2.2]{Wu88}). 
However, in the infinite dimensional case, the product of two positive operators may not be similar to a positive operator (see \cite{Ra}, \cite[Example 2.11]{Wu89}). For more development in this direction, one may consult \cite{Ra, Wu88, Wu89}.

In this paper, we study the problem when a bounded linear operator $T$ on a complex Hilbert space $H$ can be written as a product of two positive contractions. In this case, $T$ must be a contraction, and we have that
$$ -I/8\leq \text{Re} \ T  \quad {and} \quad -I/4\leq \text{Im} \ T\leq I/4$$
(see \cite[Theorem 1.1 and Corollary 4.3]{Fu}). In Proposition \ref{fact}, we give a necessary condition for this problem when $T$ has operator matrix
$$\begin{pmatrix} T_1 & T_3\\ 0 & T_2\cr\end{pmatrix}  \quad \text{on} \quad   H\oplus K.$$
In such as a case, $T_1$ and $T_2$ must also be products of two positive contractions. This is an extension of the result of Wu in \cite[Corollary 2.3]{Wu88} concerning the
finite dimensional case. 
However, even for a $2\times 2$ matrix $C$, it is not easy to determine when it is the product of two positive contractions. For example, consider
$$C = \frac{1}{25}\begin{pmatrix}9 & 3 \cr 0 & 16\cr\end{pmatrix}.$$
The diagonalizable contraction $C$ is similar to a positive operator. Thus it is a product of two 
positive operators. Moreover, $C$ satisfies $-I/8 \le {\rm Re}\,C$ and $-I/4 \le  {\rm Im}\,C 
\le I/4.$  However, we will see that $C$ cannot be written as a product of two positive 
contractions by Lemma \ref{2x2}.

Let $B(H)$ be the algebra of bounded linear operators acting on a complex Hilbert space $H$. We identify $B(H)$ with
$M_n$, the algebra of $n\times n$ complex matrices, if $H$ has finite dimension $n$. Recall that a bounded linear operator $T \in B(H)$ is positive (resp., strictly positive) if $\langle Th, h\rangle\geq 0$ (resp., $\langle Th, h\rangle>0$) for every $h\neq 0$ in $H$. We write as usual $T\geq 0$ (resp., $T>0$) when $T$ is positive (resp., strictly positive).

Recall that we call $T \in B(H)$ a quadratic operator if $(T-aI)(T-bI)=0$ for some scalars $a$, $b\in \mathbb{C}$. Every quadratic operator $T\in B(H)$ is unitarily similar to $$aI \oplus bI \oplus\begin{pmatrix} aI & P \cr 0 & bI \cr \end{pmatrix} \quad \text{on} \quad  H_1\oplus H_2\oplus (H_3\oplus H_3)$$
for some $a,b\in \mathbb{C}$, $P>0$ (see \cite{Tso}). In this paper, we  prove the following.

\begin{theorem} \label{main}
A quadratic operator $T\in B(H)$ with operator matrix
$$aI \oplus bI \oplus\begin{pmatrix} aI & P \cr 0 & bI \cr \end{pmatrix} \quad \text{on} \quad  H_1\oplus H_2\oplus (H_3\oplus H_3)$$
for some $a,b\in \mathbb{C}$ and $P>0$, can be written as a product of two positive contractions if and only if
$$a, b\in [0,1], \ \text{and} \quad  \|P\| \le  |\sqrt{a} - \sqrt{b}|\sqrt{(1-a)(1-b)}.$$
\end{theorem}

\section{Proof}

First we consider the $2\times 2$ case so that we can identify
$B(H) = M_2$ and $H = \mathbb{C}^2$.

\begin{lemma} \label{2x2}
Suppose $C = \begin{pmatrix} a & z \cr 0 & b \cr \end{pmatrix}$ with $z \ge 0$.
Then $C$ is a product of two positive contractions if and only if $a$, $b\in [0,1]$,
$$z\in S = \{c: 0 \le c \le |\sqrt{a} - \sqrt{b}| \sqrt{(1-a)(1-b)}\}.$$
If the above equivalent conditions hold, then there are continuous maps $a_{ij}(z), b_{ij}(z)$ for $1 \le i, j \le 2$ with
\begin{equation} \label{condition}
\left\{
         \begin{array}{ll} 0\leq a_{ii}(z), \ b_{ii}(z)\leq 1, \ \ a_{12}(z)=a_{21}(z)\geq 0, \ \ b_{12}(z)=b_{21}(z)\leq 0,   \\
                                              \\
         0\leq (a_{ij}(z))\leq I, \ \ 0\leq(b_{ij}(z))\leq I.
         \end{array}
         \right.
\end{equation}
such that
\begin{equation} \label{AB}
(a_{ij}(z))(b_{ij}(z)) = \begin{pmatrix}a & z \cr 0 & b \cr\end{pmatrix},  \quad
z \in S.
\end{equation}
\end{lemma}

{\em Proof}.
We first prove the sufficiency. Without loss of generality, we may assume $0\leq a\leq b\leq1$.
If  $a=b$ or $b=1$, then $z = 0$ and $C = \diag(a,1) \diag(1,b)$.
In the following, we may assume $0\leq a<b<1$, and consider two cases.

{\bf Case 1.}   $0=a<b<1$.   For $z\in S$, we have that  $z^2\leq b(1-b)$ and hence
$(z^2/b) + b \le (1-b) + b = 1$.
Consider
$$A = \begin{pmatrix} a_{11}(z) & a_{12}(z) \cr a_{21}(z) & a_{22}(z) \cr
\end{pmatrix} = \begin{pmatrix} z^2/b & z \cr z & b \cr \end{pmatrix}
\quad \hbox{ and } \quad
B = \begin{pmatrix} b_{11}(z) & b_{12}(z) \cr b_{21}(z) & b_{22}(z) \cr
\end{pmatrix} = \begin{pmatrix} 0 & 0 \cr 0 & 1 \cr \end{pmatrix},$$
Then $A$ is rank 1 with eigenvalue $(z^2/b)+b$, and $C = AB$.
Evidently, $a_{ij}(z), b_{ij}(z)$ are continuous maps for
$1 \le i, j \le 2$ and satisfy (\ref{condition}), (\ref{AB}).

{\bf Case 2.}  $0<a<b<1$.   For $z\in S$, we have
$$a+b-\frac{z^2}{(1-a)(1-b)}\geq a+b-(\sqrt{a}-\sqrt{b})^2=2\sqrt{ab}.$$
Let $\lambda_1(z) \ge \lambda_2(z)$ be
roots of the equation
$$\lambda^2-(a+b-\frac{z^2}{(1-a)(1-b)})\lambda+ab=0.$$
Then, $a\leq \lambda_2(z) \le \lambda_1(z)\leq b$ and
$\lambda_1(z), \lambda_2(z)$ are continuous maps on $z\in S$. Note that
$$\lambda_1(z)\lambda_2(z)=ab, \quad \lambda_1(z)+\lambda_2(z)=a+b-\frac{z^2}{(1-a)(1-b)}.$$
We have
\begin{equation}\label{z}
z=\sqrt{\frac{(1-a)(1-b)(\lambda_j-a)(b-\lambda_j)}{\lambda_j}}, \quad j = 1,2.\end{equation}
We will construct
$$A = \begin{pmatrix} a_{11}(z) & a_{12}(z) \cr a_{21}(z) & a_{22}(z) \cr
\end{pmatrix}= \begin{pmatrix} a_1 & a_2 \cr
a_2 & a_3\cr \end{pmatrix} \quad \hbox{ and } \quad
B = \begin{pmatrix} b_{11}(z) & b_{12}(z) \cr b_{21}(z) & b_{22}(z) \cr
\end{pmatrix}= \gamma \begin{pmatrix}  a_3 & - a_2 \cr
-a_2 & a_4\cr \end{pmatrix}$$
such that $A$ has eigenvalues $1, \lambda_1$,
$B$ has eigenvalues $1, \lambda_2$, and
$C = AB$. First, we set
\begin{equation}\label{gamma}
\gamma = \frac{\lambda_2}{b} = \frac{a}{\lambda_1} < 1.
\end{equation}
Because $1-b-\gamma + b \gamma = (1-b)(1-\gamma) > 0$, we can let
$$ a_3 = \frac{b-a}{1+b\gamma - \gamma -a}  < \frac{b-a}{b-a} = 1$$
so that by (\ref{gamma}),
\begin{eqnarray*}
a_3 - \lambda_1
&=& \frac{(b-a)}{(1+b\gamma -\gamma -a)} - \frac{a}{\gamma}
= \frac{\gamma b - \gamma a - a - \gamma ab + \gamma a + a^2}
{\gamma(1+b\gamma -\gamma - a)} \\
&=& \frac{\frac{1}{\gamma}(\gamma b - a)(1-a)}{(1+b\gamma -\gamma - a)}
= \frac{(b-\lambda_1)(1-a)}{(1+b\gamma -\gamma - a)} \geq 0.
\end{eqnarray*}
Then we can let
$$a_1 = 1+\lambda_1 - a_3 > 0  \quad \hbox{ so that } \quad a_1 + a_3 = 1+\lambda_1$$
and
$$
a_2 = \sqrt{a_1 a_3 - \lambda_1} = \sqrt{(1+\lambda_1-a_3)a_3 - \lambda_1}
= \sqrt{(1-a_3)(a_3-\lambda_1)} \quad \hbox{ so that } \quad a_1a_3-a_2^2 = \lambda_1.$$
As a result, $a_1 + a_3 = 1+\lambda_1$,  $\det(A) = \lambda_1$, and hence
$A$ has eigenvalues $1, \lambda_1$.
Further, let
$$a_4 = \frac{1}{a_3}(\frac{\lambda_2}{\gamma^2} + a_2^2)
\quad \hbox{ so that } \quad
\gamma^2(a_3a_4 - a_2^2) = \lambda_2.$$
Then by (\ref{gamma}),
\begin{eqnarray*}
\gamma(a_3+a_4)
&=&  \gamma a_3 + \frac{\gamma}{a_3}\left(\frac{\lambda_2}{\gamma^2} + a_2^2\right)
= \frac{\gamma}{a_3} \left(\frac{\lambda_2}{\gamma^2} + (a_3-\lambda_1 + \lambda_1a_3)\right)\\
&=& \frac{\gamma}{a_3}\left(\frac{\lambda_2}{\gamma^2} - \lambda_1\right) + \gamma(1+\lambda_1)
= \frac{\gamma}{a_3}\frac{(b-a)}{\gamma} +  \gamma(1+\lambda_1)\\
& = & 1 + b\gamma - \gamma - a + \gamma + \gamma \lambda_1 = 1+\lambda_2.
\end{eqnarray*}
As a result, $\tr B = 1+\lambda_2$ and $\det(B) = \lambda_2$.
Therefore, $B$ has eigenvalues $1, \lambda_2$.
Denote by $(AB)_{ij}$ the $(i,j)$ entry of $AB$. By (\ref{gamma}),
$$(AB)_{11} = \gamma(a_1a_3 - a_2^2)  = \gamma \lambda_1 = a, \quad
(AB)_{22} = \gamma(a_3a_4 - a_2^2)  = \gamma (\lambda_2/\gamma^2) = b.$$
Clearly, $(AB)_{21} = \gamma(a_2a_3-a_3a_2) = 0.$
By (\ref{gamma}) and (\ref{z}),
\begin{eqnarray*}
(AB)_{12} & = & \gamma a_2(a_4-a_1) = \gamma
\sqrt{(1-a_3)(a_3-\lambda_1)}((a_3+a_4)- (a_3+a_1)) \\
&=&  \frac{\gamma \sqrt{(1-b - \gamma + b\gamma)(b-\lambda_1)(1-a)} }
{(1+b\gamma -\gamma - a)}\left(\frac{(1+\lambda_2)}{\gamma} - (1+\lambda_1)\right)\\
&=&  \frac{ \sqrt{(1-b)(1-\gamma)(1-a)(b-\lambda_1)}}{(1+b\gamma -\gamma - \gamma\lambda_1)}
\left( {1+\lambda_2 - \gamma - \gamma \lambda_1} \right) \\
&=&  \sqrt{(1-b)(1-a)(1-\gamma)(b-\lambda_1)}
= \sqrt{ \frac{(1-b)(1-a)(\lambda_1-a)(b-\lambda_1)}{\lambda_1}} = z.
\end{eqnarray*}

%
%
%
%

For the converse, since $A$, $B$ are positive contractions with $\sigma(C)=\sigma(AB)=\sigma(B^{1/2}AB^{1/2})\subseteq [0, \infty)$, we have $0\leq a, b\leq 1$. Without loss of generality, we may assume $a\leq b$. First, consider $\|A\|=\|B\|=1$. Then the assumption $C=AB$ implies $C$ is unitarily similar to
$$\begin{pmatrix}\alpha_1 & 0 \\ 0 & 1 \end{pmatrix} \begin{pmatrix}b_1 & b_2 \\ b_2 & b_4 \end{pmatrix}=\begin{pmatrix} \alpha_1 b_1 & \alpha_1 b_2 \\ b_2 & b_4 \end{pmatrix},$$
where
$\begin{pmatrix}b_1 & b_2 \\ b_2 & b_4 \end{pmatrix}$ is unitarily similar to $\begin{pmatrix} \alpha_2 & 0 \\ 0 & 1 \end{pmatrix}$ for some $0\leq \alpha_1, \alpha_2\leq 1$, $\alpha_2\leq b_1, b_4\leq 1$ and $b_2\geq0$. Thus we have $1+\alpha_2= b_1+b_4$, $a+b=\alpha_1b_1+b_4$, $ab=\alpha_1\alpha_2=\alpha_1(b_1b_4-b_2^2)$, and $a^2+b^2+z^2=\alpha_1^2(b_1^2+b_2^2)+b_2^2+b_4^2$. These imply that
$$z^2=[\alpha_1^2(b_1^2+b_2^2)+b_2^2+b_4^2]-[(\alpha_1b_1+b_4)^2-2\alpha_1 \alpha_2]=(1-\alpha_1)^2b_2^2.$$
Hence we may assume $\alpha_1<1$. In addition, we also obtain that
$$a+b=\alpha_1b_1+b_4=\alpha_1b_1+1+\alpha_2-b_1=1+\alpha_2-(1-\alpha_1)b_1$$
and hence
\begin{eqnarray*}
b_1 & = & \frac{1}{1-\alpha_1}(1+\alpha_2-a-b)\\
& = & \frac{1}{1-\alpha_1}[(1-a)(1-b)-ab+\alpha_2]\\
& = & \frac{1}{1-\alpha_1}[\alpha_2(1-\alpha_1)+(1-a)(1-b)],
\end{eqnarray*}
where the last equality follows from $ab=\alpha_1 \alpha_2$. Let $c=(1-a)(1-b)/(1-\alpha_1)$.
Then $b_1=\alpha_2+c$ and $b_4=1-c$. By a direct computation, we see  that
\begin{eqnarray*}
z^2 & = & (1-\alpha_1)^2 b_2^2=(1-\alpha_1)^2(b_1b_4-\alpha_2)\\
& = & (1-\alpha_1)^2[(\alpha_2+c)(1-c)-\alpha_2] \qquad (\hbox{because } \alpha_2=b_1b_4-b_2^2) \\
& = & c(1-\alpha_1)[(1-\alpha_1)(1-\alpha_2)-c(1-\alpha_1)]\\
& = & (1-a)(1-b)[(a+b)-(\alpha_1+\alpha_2)],
\end{eqnarray*}
where the last equality follows from $c=(1-a)(1-b)/(1-\alpha_1)$ and $ab=\alpha_1 \alpha_2$.
Since $ab=\alpha_1 \alpha_2$,
we have $\alpha_1+\alpha_2\geq 2\sqrt{\alpha_1 \alpha_2}=2\sqrt{ab}$.
This implies that
$$z\leq |\sqrt{a} - \sqrt{b}|\sqrt{(1-a)(1-b)}.$$

In general, since $C=\alpha \begin{pmatrix} \frac{a}{\alpha} & \frac{z}{\alpha} \cr 0 &  \frac{b}{\alpha} \cr \end{pmatrix} = \alpha(\frac{A}{\|A\|})(\frac{B}{\|B\|})$, where $0<\alpha=\|A\|\|B\|\leq1$, the scalars $a, b, z$ in the above can be replaced by $a/{\alpha}, b/{\alpha}, z/{\alpha}$, respectively, to get $0\leq a/{\alpha}$, $b/{\alpha}\leq 1$ and
$$\frac{z}{\alpha}\leq \sqrt{(1-\frac{a}{\alpha})(1-\frac{b}{\alpha})}|\sqrt{\frac{a}{\alpha}}-\sqrt{\frac{b}{\alpha}}|.$$
This shows that $0\leq a$, $b\leq\alpha\leq 1$ and
$$z\leq |\sqrt{a}-\sqrt{b}|
\sqrt{(\alpha-a)(1-\frac{b}{\alpha})}\leq |\sqrt{a}-\sqrt{b}|\sqrt{(1-a)(1-b)}.$$
This proves the necessity.  {\hfill $\blacksquare$}

\medskip

In order to prove Theorem \ref{main}, we need the following fact; see, for example, \cite[p. 547]{Fo}.

\begin{lemma} \label{non}
Let $A$ be a bounded linear operator of the form
$$\begin{pmatrix}A_{11} & A_{12}\\ A_{12}^{*} & A_{22}\cr\end{pmatrix}
\ \text{on} \ H\oplus K,$$
where $H$ and $K$ are Hilbert spaces. Then $A$ is positive if and only if $A_{11}$ and $A_{22}$
are both positive and there exists a contraction $D$ mapping $K$ into $H$ satisfying
$A_{12}=A_{11}^{1/2}DA_{22}^{1/2}$.
\end{lemma}

\begin{lemma} \label{spectral}
Suppose $a_{11}(z), a_{22}(z), a_{12}(z) = a_{21}(z)$ 
are continuous real-valued functions defined on $S \subseteq [0, \infty)$ such that
$A=\begin{pmatrix} a_{11}(z) &  a_{12}(z) \cr  a_{21}(z) & a_{22}(z) \cr \end{pmatrix} 
\geq
0$ for all $z\in S$.
Then $\begin{pmatrix} a_{11}(P) & a_{12}(P) \cr a_{21}(P) & a_{22}(P) \cr \end{pmatrix}\geq 0$
on $H\oplus H$ for all positive operators $P\in B(H)$ with spectrum in $S$.
\end{lemma}

{\em Proof}.
Since $A\ge 0$, we have $a_{11}(z), a_{22}(z)\ge 0$ and
$$0\le a_{12}(z)a_{21}(z)\leq a_{11}(z)a_{22}(z), \quad   \quad  z\in S.$$
Define $h(z)$ by
$$h(z):= \left\{
         \begin{array}{ll}
          \frac{a_{12}(z)}{a_{11}^{1/2}(z)a_{22}^{1/2}(z)}& \mbox{if }  |a_{12}(z)| >0,  \\
         0  & \mbox{if }  a_{12}(z) =0.
         \end{array}
         \right.$$
Then $h(z)$ is a bounded Borel function on $S$ with $|h(z)|\leq 1$, which satisfies $$a_{12}(z)=a_{11}^{1/2}(z)h(z)a_{22}^{1/2}(z).$$
By the spectral theorem,
for all positive operators $P\in B(H)$ with spectrum in $S$, 
we have $a_{11}(P)\geq 0$, $a_{22}(P)\geq 0$, $a_{12}(P)=a_{21}(P)\geq 0$ and
$$a_{12}(P)=a_{11}^{1/2}(P)h(P)a_{22}^{1/2}(P)$$
for the contraction $h(P)\in B(H)$. 
Our assertion follows from Lemma \ref{non}.   {\hfill $\blacksquare$}

%
%

In the finite dimensional case, Wu \cite[Corollary 2.3]{Wu88} has shown that if $C=\begin{pmatrix} C_1 & C_3\\ 0 & C_2\cr\end{pmatrix}$ is a product of two positive operators, then so are $C_1$ and $C_2$. Proposition \ref{fact} gives another proof which holds for both finite and infinite dimensional Hilbert spaces. In fact, it is also true that positive operators are replaced by positive contractions.

\begin{proposition} \label{fact}
Let $T$ be a bounded linear operator of the form
$$\begin{pmatrix} T_1 & T_3\\ 0 & T_2\cr\end{pmatrix} \ \text{on} \ H\oplus K,$$
where $H$ and $K$ are both Hilbert spaces. If $T$ is a product of
two positive contractions, then so are $T_1$ and $T_2$.
\end{proposition}


{\em Proof}.
By our assumption and Lemma \ref{non}, we may assume that $T=AB$, where $A$ and $B$ are
of the form
$$\begin{pmatrix} A_1 & A_1^{1/2}D_1A_2^{1/2}\cr
A_2^{1/2}D_1^{*}A_1^{1/2} & A_2 \cr \end{pmatrix} \quad
\hbox{ and } \quad
\begin{pmatrix} B_1 & B_1^{1/2}D_2B_2^{1/2} \cr
B_2^{1/2}D_2^{*}B_1^{1/2} & B_{2} \cr\end{pmatrix}
\quad \text{on} \ H\oplus K,$$
respectively, such that  $0\leq A_1\leq I_H$, $0\leq A_2\leq I_K$,
$0\leq B_1\leq I_H$, $0\leq B_2\leq I_K$,
$D_1$ and  $D_2$ are contractions from $K$ into $H$. From $T=AB$, we obtain that
\begin{eqnarray}\label{eq:1}
T_1=A_1B_1+A_1^{1/2}D_1(A_2^{1/2}B_2^{1/2}D_2^{*}B_1^{1/2}),
\end{eqnarray}
\begin{eqnarray}\label{eq:2}
A_2^{1/2}(D_1^{*}A_1^{1/2}B_1^{1/2})B_1^{1/2}=-A_2^{1/2}(A_2^{1/2}B_2^{1/2}D_2^{*})B_1^{1/2},
\end{eqnarray}
\begin{eqnarray*}
T_2=(A_2^{1/2}D_1^{*}A_1^{1/2}B_1^{1/2})D_2B_2^{1/2}+A_2B_2.
\end{eqnarray*}
Let $E_1$ be the restriction of $A_2^{1/2}$ to (\ker $A_2^{1/2})^{\perp}$, then $E_1$ is injective. Since $0\leq A_2^{1/2}\leq I_K$, so we can consider the (possibly unbounded) inverse $E:=E_1^{-1}$: ran $A_2^{1/2}$ $\rightarrow$ (\ker $A_2^{1/2})^{\perp}$ such that $EA_2^{1/2}=P_{\overline{ran A_2^{1/2}}}$. Hence by (\ref{eq:2}), we derive that
$$A_2^{1/2}B_2^{1/2}D_2^{*}B_1^{1/2}=P_{\overline{ran A_2^{1/2}}}(A_2^{1/2}B_2^{1/2}D_2^{*}B_1^{1/2})=-P_{\overline{ran A_2^{1/2}}}(D_1^{*}A_1^{1/2}B_1).$$
Moreover, substitute this into (\ref{eq:1}) to get
\begin{eqnarray*}
T_1 & = & A_1B_1-A_1^{1/2}D_1(P_{\overline{ran A_2^{1/2}}}(D_1^{*}A_1^{1/2}B_1))\\
& = & [A_1^{1/2}(I_H-D_1P_{\overline{ran A_2^{1/2}}}D_1^{*})A_1^{1/2}]B_1\\
& = & [A_1^{1/2}(I_H-(P_{\overline{ran A_2^{1/2}}}D_1^{*})^{*}
(P_{\overline{ran A_2^{1/2}}}D_1^{*}))A_1^{1/2}]B_1.
\end{eqnarray*}
Note that $\|P_{\overline{ran A_2^{1/2}}}D_1^{*}\|\leq 1$ implies that
$$0\leq (I_H-(P_{\overline{ran A_2^{1/2}}}D_1^{*})^{*}(P_{\overline{ran A_2^{1/2}}}D_1^{*}))
\leq I_H.$$
Therefore, $T_1=[(A_1^{1/2}P_1^{*})P_1A_1^{1/2}]B_1$,
where $P_1^{*}P_1=I_H-(P_{\overline{ran A_2^{1/2}}}D_1^{*})^{*}
(P_{\overline{ran A_2^{1/2}}}D_1^{*})$ for some positive contraction $P_1$ on $H$.
This shows that $T_1$ is a product of two positive contractions. Similarly, we can show that
$T_2^{*}$ is a product of two positive contractions, and hence so is $T_2$. This
completes our proof.    {\hfill $\blacksquare$}

\medskip

Now we are ready to give the proof of Theorem \ref{main}.

{\em Proof of Theorem \ref{main}}. 
We first prove the necessity. By assumption, we can focus on the part $$\begin{pmatrix} aI & P \cr 0 & bI \cr \end{pmatrix}\in B(H_3\oplus H_3)$$ for some $P>0$ . Now, consider a $2\times 2$ matrix $\begin{pmatrix} a & z \cr 0 & b \cr \end{pmatrix}$ with $a$, $b\in [0,1]$ and
$$z\in S := \{c:
0 \le c \le |\sqrt{a} - \sqrt{b}| \sqrt{(1-a)(1-b)}\}.$$
Then by Lemma \ref{2x2}, there are continuous maps $a_{ij}(z), b_{ij}(z)$ for $1 \le i, j \le 2$ with $a_{12}(z)=a_{21}(z)\ge 0$, $b_{12}(z)=b_{21}(z)\le 0$ and satisfy
$$0\leq (a_{ij}(z))\le I_2, \ \ 0\leq (b_{ij}(z))\le I_2,    \quad  \quad  (a_{ij}(z))(b_{ij}(z)) = \begin{pmatrix}a & z \cr 0 & b \cr\end{pmatrix},     \quad z\in S.
$$
By Lemma \ref{spectral},
$$ 0 \le (a_{ij}(P)) \le I
\quad  \hbox{ and } \quad
0 \le (b_{ij}(P)) \le I.$$
By the spectral theorem on positive operators,
$$(a_{ij}(P))(b_{ij}(P)) = \begin{pmatrix}aI & P \cr 0 & bI \cr\end{pmatrix}.$$


%
%
%
%
%
%

To prove the converse, suppose there is a factorization of the quadratic operator $T\in B(H)$
with operator matrix $aI \oplus bI \oplus \begin{pmatrix} aI & P \cr 0 & bI \cr \end{pmatrix}$
for some $P \ge 0$ as the product
of two positive contractions. By Proposition \ref{fact}, we know that
$$T_1=\begin{pmatrix} aI & P \cr 0 & bI \cr \end{pmatrix}=AB \quad \quad \quad \quad
\text{for \ some} \quad 0 \leq A, B\leq I, \ A, B\in B(H_3\oplus H_3).$$
We may use the Berberian construction (see \cite{Be})
to embed $H_3$ into a larger Hilbert space $K_3$, $B(H_3)$ into $B(K_3)$.
Suppose $A = (A_{ij})_{1 \le i, j \le 2}, B = (B_{ij})_{1 \le i, j \le 2} \in B(H_3\oplus H_3)$.
Then $P$, $A$, and $B$ are extended to $\tilde P \in B(K_3)$,
$\tilde A = (\tilde A_{ij})_{1 \le i, j \le 2}\in B(K_3\oplus K_3)$, and 
$\tilde B = (\tilde B_{ij})_{1 \le i, j \le 2} \in B(K_3\oplus K_3)$, respectively,
such that the following conditions hold.

\medskip
(a)  $\tilde P\ge0$ with $\|P\| = \|\tilde P\|$ such that all the
elements in $\sigma(\tilde P)$ are eigenvalues of $\tilde P$.

(b)   $0 \leq \tilde A, \tilde B\leq I$ such that 
$\tilde T_1=\begin{pmatrix} aI & \tilde P \cr 0 & bI \cr \end{pmatrix}=\tilde A\tilde B$.

\medskip\noindent
Since $\tilde P \ge 0$ and $\sigma(\tilde P)$ are eigenvalues of $\tilde P$,
the quadratic operator  $\tilde T_1$ is unitarily similar to
$\begin{pmatrix} a & \|P\| \cr 0 & b \cr \end{pmatrix} \oplus T_2$ that admits a factorization
as the product of two positive contractions.
By Proposition \ref{fact}, we see that
$\begin{pmatrix} a & \|P\| \cr 0 & b \cr \end{pmatrix}$
is a product of two positive contractions. Thus,
$$\|P\| \le  |\sqrt{a} - \sqrt{b}|\sqrt{(1-a)(1-b)}.$$
\vskip -.3in {\hfill $\blacksquare$}

  \medskip

\section*{acknowledgment}
Li is an honorary professor of the University of Hong Kong and the Shanghai University. His research was supported by US NSF and HK RCG. The Research of Tsai was supported by the National Science Council of the Republic of China under the project NSC 102-2811-M-110-018. Tsai would like to thank Pei Yuan Wu and Ngai-Ching Wong  for their helpful suggestions and comments. Some results in this paper are contained in the doctorial thesis of Ming-Cheng Tsai under the supervisor of Pei Yuan Wu to whom Tsai would express his heartfelt thanks.

\end{document}